\input amstex
\input amssym
\documentstyle{amsppt}
\magnification=1200%


\def\a{\alpha}
\def\b{\beta}

\def\d{\delta}

\def\s{\sigma}

\def\x{\times}

\def\o{\overline}
\def\f{\flushpar}

\def\B{\Cal B}

\def\T{\widehat T}
\def\({\biggl(}
\def\){\biggr)}\def\cempt{\text{conservative, ergodic, measure preserving transformation}}

\def\mpt{\text{measure preserving transformation}}
\def\<{\bold\langle}
\def\>{\bold\rangle}

\def\bul{\smallskip\f$\bullet\ \ \ $}\def\sms{\smallskip\f}\document
\topmatter\title Absolutely continuous, invariant measures
 for dissipative, ergodic transformations\endtitle\author Jon. Aaronson,\   Tom Meyerovitch\endauthor
 \address[Jon. Aaronson]{\ \ School of Math. Sciences, Tel Aviv University,
69978 Tel Aviv, Israel.}
\endaddress
\email{aaro\@tau.ac.il}\endemail\address[Tom Meyerovitch]{\ \
School of Math. Sciences, Tel Aviv University, 69978 Tel Aviv,
Israel.}
\endaddress
\email{tomm\@tau.ac.il}\endemail\abstract We show that a
dissipative, ergodic $\mpt$ of a $\s$-finite, non-atomic measure
space always has many non-proportional, absolutely continuous,
invariant measures and  is ergodic with respect to each one of
these.
\endabstract\thanks\copyright 2005. Preliminary version.
\endthanks\subjclass 37A05,\ 37A40\endsubjclass
\endtopmatter\rightheadtext{invariant measures}
\heading \S0 Introduction\endheading
 Let $(X,\B,m,T)$ be an invertible, ergodic  $\mpt$ of a  $\s$-finite measure space, then there are no other
 $\s$-finite, $m$-absolutely continuous, $T$-invariant measure other than constant multiples of $m$,
 because the density of any such measure is $T$-invariant, whence constant by ergodicity.

 When $T$ is not invertible, the situation becomes more
 complicated.

If $(X,\B,m,T)$ is a  $\cempt$  of a  $\s$-finite measure space,
then (again) there are no other
 $\s$-finite, $m$-absolutely continuous, $T$-invariant measure other than constant
 multiples of $m$ (see e.g. theorem 1.5.6 in
\cite{A}). When $T$ is not conservative, the situation is different.

In this note, we show (proposition 1) that a dissipative \f$\mpt$
has many non-proportional, $\s$-finite, absolutely continuous,
invariant measures. \par If the dissipative $\mpt$ is ergodic
(exact), then it  is also ergodic  (exact) with respect to each of
these $\s$-finite, absolutely continuous, invariant measures
(proposition 2).

Proposition 1 was known for certain examples: the ``Engel series
transformation" (see\cite{T},  also \cite{S1}); the one sided shift
of a random walk on a polycyclic group with centered, adapted jump
distribution (ergodicity follows from \cite{K}, existence of
non-proportional invariant densities follows from \cite{B-E}); and
the Euclidean algorithm transformation (see \cite{D-N} which
inspired this note). More details are given in \S2.

\S1 is devoted to results (statements and proofs) and \S2 has
examples of ergodic, dissipative $\mpt$s.

To conclude this introduction, we consider  \subheading{An
illustrative example}

Fix $q\in (0,1)$ and consider the stochastic matrix $p:\Bbb Z\x\Bbb
Z\to [0,1]$ defined by $p_{s,s}:=1-q,\ p_{s,s+1}:=q$ and $p_{s,t}=0\
\forall\ t\ne s,s+1$. Let $(X,\B,m,T)$ be the {\tt one-sided Markov
shift} with $X:=\Bbb Z^\Bbb N$, $\B$ the $\s$-algebra generated by
{\it cylinders} (sets of form $[a_1,\dots,a_k]:=\{x\in X:\ x_j=a_j\
\forall \ 1\le j\le k\}$ and $m:\B\to [0,\infty]$ the measure
satisfying $m([a_1,\dots,a_k]):=\prod_{j=1}^{k-1}p_{a_j,a_{j+1}}$.
It is not hard to check that $(X,\B,m,T)$ is a $\mpt$. By random
walk theory (see \S2 and \cite{D-L}) it is {\it exact} in the sense
that $\bigcap_{n\ge 0}T^{-n}\B\overset{m}\to{=}\{\emptyset,X\}$. It
can be checked directly that $F:X\to [0,\infty)$ defined by
$$F(x_1,x_2,\dots)=\cases & 0\ \ \ \ \
N_0(x):=\sum_{n=1}^\infty\d_{x_n,0}>1,\\ & 1\ \ \ \ \ N_0(x)=1,\
x_1<0,\\ & q\ \ \ \ \ \ \ \ \text{else}\endcases$$ is the density
of a $\s$-finite, $m$-absolutely continuous, $T$-invariant
measure.

\heading\S1 Results\endheading \subheading{Wandering sets}

For a $\mpt$ $(X,\B,m,T)$ let
 $\Cal W_T:=\{W\in\B:\ W\cap T^{-n}W=\emptyset\ \forall\ n\ge 1\}$:  the
collection of {\it wandering sets} for $T$. As is well known (see
e.g. \cite{A}, or \cite{Kr}) $T$ is dissipative iff $X$ is a
countable union of wandering sets $\mod m$.

In case $T$ is dissipative, invertible then \bul$\exists\
W_{\max}\in\Cal W_{T},\ \biguplus_{n\in\Bbb Z} T^nW_{\max}= X\ \mod\
m$ (see e.g. \cite{A}, or \cite{Kr}); and \bul if $W\in\Cal W_{ T},$
then
  $\biguplus_{n\in\Bbb Z}
T^nW= X\ \mod\  m$ only if $ m(W)= m(W_{\max})$, the reverse
implication holding when $ m(W_{\max})<\infty$ (see theorem 1 in
\cite{H-K} ). We denote the constant $ m(W_{\max})=:\goth w(T)$.

\

We prove: \proclaim{Proposition 1}

 Let
$(X,\B,m,T)$ be a dissipative $\mpt$ of a standard, non-atomic,
$\s$-finite measure space, then  $\exists\ c\in (0,\infty]$ so that
for every $W\in\Cal W_T,\ m(W)<c,\ \exists$ a non-zero,
$m$-absolutely continuous, $T$-invariant measure $\mu$ with bounded
density so that $\mu(W)=0$.
\endproclaim
\demo{Proof} By Rokhlin's theorem (see \cite {Ro} or theorem 3.1.5
in \cite{A}), there is an invertible, $\mpt$ $(\widetilde
X,\widetilde\B,\widetilde m,\widetilde T)$ equipped with a
measurable map $\pi:\widetilde X\to X$ satisfying $$ \pi\circ
\widetilde T=T\circ\pi,\ \ \widetilde m\circ\pi^{-1}=m.\tag\ddag$$
It follows that  $(\widetilde X,\widetilde\B,\widetilde
m,\widetilde T)$ is dissipative. \par Given $p\in
L^\infty(\widetilde X)_+,\ p\circ\widetilde T=p$, define
$\mu_p\in\goth M(X,\B)$ by
$$\mu_p(A):=\int_X1_{A}\circ\pi p\ d\widetilde m.$$ Evidently $\mu_p\ll m$ with
$\|\tfrac{d\mu_p}{dm}\|_\infty\le \|p\|_\infty$ and
$\mu_p(T^{-1}A) =\mu_p(A)\ \ \ (A\in\B)$.

\

Next we show, as advertised, that    each wandering set of small
enough measure is annihilated by some $\mu_p$.

\

\

Let $c:=\goth w(\widetilde T)\in (0,\infty]$ and suppose that
$W\in\Cal W_T,\ m(W)<c$, then $\pi^{-1}W\in\Cal W_{\widetilde T}$
and $\widetilde m(\widetilde X\setminus\biguplus_{n\in\Bbb
Z}\widetilde T^n\pi^{-1}W)>0$.

Set $Y:=\widetilde X\setminus\biguplus_{n\in\Bbb Z}\widetilde
T^n\pi^{-1}W$ (then $\widetilde TY=Y$) and let $\mu:=\mu_{1_Y}$,
then (as above) $\mu\ll m$ with $\|\tfrac{d\mu}{dm}\|_\infty\le 1$
and $\mu(T^{-1}A) =\mu(A)\ \ \ (A\in\B)$.

By construction, $\mu(W)=\widetilde m(\pi^{-1}W\cap Y)=0.\qed$
\enddemo

 \subheading{Remarks}

 1)  The  density $F$ in the illustrative example above,  can be obtained as in the proof of proposition
 1 as
 $$\int_AFdm=\widetilde m(\pi^{-1}A\cap \biguplus_{n\in\Bbb Z}\widetilde
 T^n\pi^{-1}[-1,0,1])$$ or
 $$F=\sum_{n\ge 0}1_{[-1,0,1]}\circ T^n+\sum_{n\ge
 1}\T_m^n1_{[-1,0,1]}$$ where $\T_m$ denotes the {\it transfer operator} of the  $\mpt$
 $(X,\B,m,T)$, which  is the
operator  defined on the space $L(X)_+$ of non-negative,
measurable functions by $\int_A\T_m fdm=\int_{T^{-1}A}f dm \ \
(f\in L(X)_+,\ A\in\B)$.

2) Evidently, $p\in L(X)_+$ is the density of an $m$-absolutely
continuous, $T$-invariant measure iff $\T_m p=p$. Also $\T_m
(f\circ T)=f$.

If $(X,\B,m,T)$ is a dissipative $\mpt$, then
$$\sum_{n\ge 0}f\circ T^n<\infty\ \&\ \sum_{n\ge 0}\T_m^nf<\infty\
\forall\ f\in L^1(X),\ f\ge 0.$$ It follows that $\T_m F=F$ where
$F=F(f):=\sum_{n\ge 0}f\circ T^n+\sum_{n\ge 1}\T_m^nf$ whenever
$f\in L^1$. This can be used to prove a less precise version of
proposition 1 without assuming standardness of $(X,\B,m)$: if $A,\
B\in\B$ are disjoint and $A\uplus B\in\Cal W_T$, then
$F(1_A)1_B=0\ \mod m$.

 3) Let $(X,\B,m,T)$ be a dissipative $\mpt$ of a standard,
non-atomic, $\s$-finite measure space and let $(\widetilde
X,\widetilde\B,\widetilde m,\widetilde T)$ be its  {\it natural
extension}: i.e.   an invertible, $\mpt$ $(\widetilde
X,\widetilde\B,\widetilde m,\widetilde T)$ equipped with a
measurable map $\pi:\widetilde X\to X$ satisfying (\ddag) and a
minimality condition that
$$\bigvee_{n=1}^\infty\widetilde T^{n} \pi^{-1}\B=\widetilde \B\ \mod\ \widetilde m.$$
Natural extensions are unique up to isomorphism, and exist by
Rokhlin's theorem (mentioned above). We claim that any
$m$-absolutely continuous, $T$-invariant measure $\mu$ with
bounded density is of form $\mu_p$ where $p\in L^\infty(\widetilde
X),\ p\circ\widetilde T=p$.

To see this, let $\mu:\B\to [0,\infty]$ be such a measure. Now
define the $\widetilde T$-invariant measure $\widetilde\mu$ on
$(\widetilde X,\widetilde\B)$ as in the proof of theorem 3.1.5 in
\cite{A}. Evidently $\widetilde\mu\ll\widetilde m$,
$p:=\tfrac{d\widetilde\mu}{d\widetilde m}$ is a bounded,
measurable, $\widetilde T$-invariant function and $\mu=\mu_p$.

 \proclaim{Proposition 2}

 Let $(X,\B,m,T)$ be an ergodic (exact)
$\mpt$ of a standard,  $\s$-finite measure space.

 If
$\mu\ll m$ is a $\s$-finite, $T$-invariant measure, then
$(X,\B,\mu,T)$ also an ergodic (exact)  $\mpt$.
\endproclaim
\subheading{Remark} Proposition 2 applies mainly to dissipative,
ergodic (exact)
\f$\mpt$s of  standard, non-atomic $\s$-finite measure \f spaces.
\demo{Proof} By theorem 2 in \cite{D}, $(X,\B,\mu,T)$ is \bul
ergodic iff $\|\tfrac1n\sum_{k=0}^{n-1}\T_\mu^ku\|_{L^1(\mu)}\to
0$ for each $u\in L^1(\mu)_0$; \bul exact  iff
$\|\T_\mu^nu\|_{L^1(\mu)}\to 0$ for each $u\in L^1(\mu)_0$.

Here $L^1(\mu)_0:=\{u\in L^1(\mu):\ \int_Xud\mu=0\}$.

Suppose that $p\in L(X)_+,\ \T_mp=p$. We'll show that $(X,\B,m,T)$
exact implies that $(X,\B,\mu,T)$ is also  exact where $d\mu=pdm$.
The proof for ergodicity is analogous. We note first that
$$\T_{\mu} f=1_{[p>0]}\tfrac1p\T_m(fp).$$
Suppose that $u\in L^1(\mu)_0$, then $up\in L^1(m)_0$ and by
exactness of $(X,\B,m,T)$, $\|\T_m^n(up)\|_{L^1(m)}\to 0$. Thus
$$\|\T_\mu^nu\|_{L^1(\mu)}=\int_X1_{[p>0]}|\T^n_m(up)|dm\le\|T_m^n(up)\|_{L^1(m)}\to
0$$ and $(X,\B,\mu,T)$ is  exact.\hfill\qed\enddemo \heading\S2
Examples of ergodic, dissipative $\mpt$s\endheading
\subheading{The Engel series transformation}\par This is the
piecewise linear map $T:(0,1]\to (0,1]$ defined  by
$T(x):=([\tfrac1x]+1)x-1$ considered with respect to Lebesgue
measure. Dissipation follows from $T^nx\downarrow 0$ for each
$x\in (0,1)\setminus\Bbb Q$, ergodicity was shown in \cite{S2} and
invariant densities were given explicitly in \cite{T}. This
material is also in the book \cite{S1}. \subheading{Dissipative,
ergodic, random walks }
\par The (left) {\it random walk on LCP group $\Bbb G$} with
{\it jump probability } $p\in\Cal P(\Bbb G)$ ($\text{\tt RW}(\Bbb
G,p)$) is $(X,\B,\mu,T)$ the stationary, one-sided shift of the
Markov chain on $\Bbb G$ with transition probability
$P(g,A):=p(Ag^{-1})\ \ (A\in\B(\Bbb G))$ defined by
$$X:=\Bbb G^\Bbb N,\ \B:=\B(X),\
T(x_1,x_2,\dots)=(x_2,x_3,\dots)$$ and
$$\mu([A_1,A_2,\dots,A_N]):=\int_\Bbb
GP_x([A_1,A_2,\dots,A_N])dm(x)$$ where $m$ is a left Haar measure
on $\Bbb G$ and for $ A_1,A_2,\dots,A_N\in\B(\Bbb G)$,
$$[A_1,A_2,\dots,A_N]:=\{x=(x_1,x_2,\dots)\in X:\ x_k\in A_k\
\forall\ 1\le k\le N\};$$
$$P_x([A_1]):=1_{A_1}(x),\ P_x([A_1,A_2,\dots,A_N]):=1_{A_1}(x)\int_\Bbb
GP_{gx}([A_2,\dots,A_N])dp(g).$$
\par For an Abelian group $\Bbb G$ it is shown in \cite{D-L} (using
\cite{F}) that $\text{\tt RW}(\Bbb G,p)$ is ergodic iff
$\o{\<\text{\rm spt.}\,p\>}=\Bbb G$, and exact iff $\o{\<\text{\rm
spt.}\,p\ -\ \text{\rm spt.}\,p\>}=\Bbb G$. An exact random walk
on  $\Bbb Z^d$ can be conservative or dissipative when $d=1,2$ but
is always dissipative when $d\ge 3$.

\subheading{Dissipative, exact inner functions}
\par By Herglotz's theorem, any analytic endomorphism  $F:\Bbb R^{2+}:=\{x+iy\in\Bbb C:\
y>0\}\to\Bbb R^{2+}$ has the form
$$F(z)=\a z+\b+\int_{\Bbb R}\({1+tz\over t-z}\)d\mu(t)\tag2$$
where $\a\ge 0,\ \b\in\Bbb R$ and $\mu$ is a positive measure on
$\Bbb R$. The limits $\lim_{y\to 0+}F(x+iy)$ exist for a.e.
$x\in\Bbb R$. The analytic endomorphism $F:\Bbb R^{2+}\to\Bbb
R^{2+}$ is called an {\it inner function} if $T(x):=\lim_{y\to
0+}F(x+iy)\in\Bbb R$ for a.e. $x\in\Bbb R$, equivalently: $\mu$ is
a singular measure on $\Bbb R$. A (referenced) discussion of inner
functions can be found in chapter 6 of \cite{A}.

It is known that the real restriction $T$ of an inner function is
Lebesgue non-singular: $m(T^{-1}A)=0\ \iff\ m(A)=0\ \ (A\in\B(\Bbb
R))$ where $m$ is Lebesgue measure on $\Bbb R$ (see e.g.
proposition 6.2.2  in \cite{A}) and that $m\circ T^{-1}=m$ when
$\a=1$ in (2) (see e.g. proposition 6.2.4  in \cite{A}). If $\b=0$
and  $\mu$ is  symmetric measure ($\mu(-A)=\mu(A)$), then the real
restriction $T$ is odd, and exact by theorem 6.4.5 in \cite{A}.
\par If, in addition, $\mu([-x,x]^c)\propto\tfrac1{x^\a}$ for some
$0<\a<1$, then by lemma 6.4.7 in \cite{A}, $T$ is dissipative.

\subheading{Dissipative, ergodic, number theoretical
transformations}
\par The {\it Euclidean algorithm} is the transformation defined
by $T:\Bbb R_+^2\to\Bbb R_+^2$ by $$T(x,y)=\cases & (x-y,y)\ \ \ \
\ x>y,\\ & (x,y-x)\ \ \ \ \ x<y.\endcases$$ It is shown in
\cite{D-N} that $(\Bbb R_+^2,\B(\Bbb R_+^2),\mu,T)$ is an ergodic,
dissipative,\f  $\mpt$ where $d\mu(x,y)=\tfrac{dxdy}{xy}$.
Exactness does not seem to be known. \par The {\it Rauzy
induction} transformations considered in \cite{V} are also known
to be ergodic, dissipative $\mpt$s.
\subheading{ Dissipative $S$-unimodal maps} \sms These are discussed in [B-H] in terms of their attractors. Conditions are given for ergodicity, exactness, dissipativity and existence of  $\s$-finite
invariant densities.

 \heading References\endheading

\enddocument